\documentstyle[amssymb,11pt,fleqn]{article}
\def\la{\langle}\def\ra{\rangle}

\def\pf{{\bf Proof\quad }}
\def\pfend{\hfill{$\Box$}\vskip 0.2cm}

\begin{document}

\title{A generalization of $c$-Supplementation
\thanks{Supported by the National Natural Science
Foundation of China (Grant No. 10171074).}
\author{Shiheng Li, Dengfeng Liang, Wujie Shi\\
\small School of Mathematic Science, Suzhou University, Suzhou 215006, China\\
\small E-mail: lishiheng01@163.com; dengfengliang@163.com; wjshi@suda.edu.cn}}

\date{}
\maketitle

\begin{abstract} A subgroup $H$ is said to be $nc$-supplemented
in a group $G$ if there is a subgroup $K\leq G$ such that
$HK\unlhd G$ and $H\cap K$ is contained in $H_G$, the core of $H$
in $G$. We characterize the solvability of finite groups $G$ with
some subgroups of Sylow subgroups $nc$-supplemented in $G$. We
also give a result on $c$-supplemented subgroups.
\end{abstract}

\textbf{Keywords}\,\, solvable, $nc$-supplemented, $p$-nilpotent.

Wang\cite{FC} introduces the notation of $c$-supplemented
subgroups and determines the structure of finite groups $G$ with
some subgroups of Sylow subgroups $c$-supplemented in $G$. Here we
give a new concept called $nc$-supplementation that is a weak
version of $c$-supplementation and characterize the solvability of
groups $G$ with some maximal or 2-maximal subgroups
$nc$-supplemented in $G$, respectively.

In this paper, $\pi$ denotes a set of primes. We say $G\in E_\pi$
if $G$ has a Hall $\pi$-subgroup; $G\in C_\pi$ if $G\in E_\pi$ and
any two Hall $\pi$-subgroups of $G$ are conjugate in $G$; $G\in
D_\pi$ if $G\in C_\pi$ and every $\pi$-subgroup of $G$ is
contained in a Hall $\pi$-subgroups of $G$. We say that a number
$n$ is a $\pi$-number if every of its prime divisor is in $\pi$.
$|G|_\pi$ denotes the largest $\pi$-number that divides $|G|$.
$H<\cdot G$ denotes that $H$ is a maximal subgroup of $G$. $L$ is
called a 2-maximal subgroup of $G$, if there exists a maximal
subgroup $M$ of $G$ such that $L<\cdot M$.\\

{\bf Definition 1.}\ Let $G$ be a group and $H$ be a subgroup of
$G$.

(1) $H$ is said to be $c$-supplemented in $G$ if there exists a
subgroup $K$ of $G$ such that $HK=G$ and $H\cap K\leq H_G$, where
$H_G=\bigcap\limits_{g\in G}H^g$ is the core of $H$ in $G$. We say
that $K$ is a $c$-supplement of $H$ in $G$.

(2) $H$ is said to be $nc$-supplemented in $G$ if there is a
subgroup $K$ of $G$ such that $HK\unlhd G$ and $H\cap K\leq H_G$.
We say that $K$ is a $nc$-supplement of $H$ in $G$.\pfend

{\bf Remark 1}. Let $H$ be $nc$-supplemented in $G$. Evidently,
$H$ is $c$-supplemented in $G$ if $G$ is simple; $H$ is
$c$-supplemented in $G$ if $H$ is a maximal subgroup of $G$.\pfend

In general, $nc$-supplementation does not implies
$c$-supplementation.

{\bf Example.} Let $G=A_4$ and
$B=\{(1),(12)(34),(13)(24),(14)(23)\}$. Set $C=\{(1),(12)(34)\}$
and $D=\{(1),(14)(23)\}$. Then $C\times D=B\unlhd G$ and $C$ is
$nc$-supplemented in $G$. However, $C$ is not $c$-supplemented
since $C_G=1$ and $A_4$ has no subgroup of order 6.\pfend

{\bf Lemma 1.}\ If $H$ is $nc$-supplemented in $G$, then there
exists a subgroup $C$ of $G$ such that $H\cap C=H_G$ and $HC\unlhd
G$.

\pf Suppose $H$ is $nc$-supplemented in $G$. Then there is a
subgroup $C_1\leq G$ such that $H\cap C_1\leq H_G$ and $HC_1\unlhd
G$. Let $C=C_1H_G$. Then $HC=(HC_1)H_G\unlhd G$ and $H\cap C=H\cap
C_1H_G=H_G(H\cap C_1)=H_G$.\pfend

{\bf Lemma 2.} Let $H$ be $nc$-supplemented in $G$.

(1) If $H\leq H\leq G$, then $H$ is $nc$-supplemented in $M$.

(2) If $N\unlhd G$ and $N\leq H$, then $H/N$ is $nc$-supplemented
in $G/N$.

(3) If $N\unlhd G$ and $(|N|,|H|)=1$, then $HN/N$ is
$nc$-supplemented in $G/N$.

\pf Similar to the argument in the proof of \cite[Lemma 2.1 (1)
and (2)]{FC}, we get (1) and (2), respectively.

(3) If $H$ is $nc$-supplemented in $G$, then there is a subgroup
$C\leq G$ such that $HC\unlhd G$ and $H\cap C\leq H_G$. Evidently,
$(HN/N)(CN/N)=(HC)N/N\unlhd G/N$. On the other hand, since
$(HN\cap C)H=HN\cap HC\leq HN$, $(HN\cap C)N\leq HN$ and
$(|N|,|H|)=1$, we get $HN\cap C=(H\cap C)(N\cap C)$ by \cite[A,1.6
Lemma,(c)]{DH}. Hence $HN/N\cap CN/N=(H\cap C)N/N\leq H_GN/N\leq
(HN/N)_{G/N}$. This gives that $HN/N$ is $nc$-supplemented in
$G/N$.\pfend

{\bf Lemma 3}\cite[Proposition 2.1]{AF}. If $K$ is a normal
subgroup of the group $G$ such that $K\in C_\pi$ and $G/K\in
E_\pi$, then $G\in E_\pi$.\pfend

From \cite[p.485, Theorem]{FG} we get the following result by
computing $|G|_2$:

{\bf Lemma 4.} Let $G$ be a simple group having a Sylow 2-subgroup
isomorphic to $C_2\times C_2$. Then $G\cong L_2(q)$, where
$q\equiv 3\ (mod\ 8)$ or $q\equiv 5\ (mod\ 8)$.\pfend

{\bf Theorem 1.} Let $G$ be a finite group and let $P$ be a Sylow
$p$-subgroup of $G$, where $p$ is a prime divisor of $|G|$.
Suppose that there is a maximal subgroup $P_1$ of $P$ such that
$P_1$ is $nc$-supplemented in $G$.

(a) If $P_1\neq 1$, then $G$ is not a non-Abelian simple group.

(b) If $p=2$, then $G\in E_{2'}$ and every composition factor of
$G$ is either a cyclic group of prime order or isomorphic to
$L_2(r)$, where $r=2^n-1$ is a Mersenne prime.

\pf (a) Suppose $P_1$ is $nc$-supplemented in $G$. If $G$ is
simple then $P_1$ is $c$-supplemented in $G$. By \cite[Theorem
2.2]{FC}, it follows that $G$ is not simple, a contradiction.

(b) Suppose $p=2$. If $(P_1)_G\neq 1$, then $G/(P_1)_G$ satisfies
the hypothesis by Lemma 2. Hence $G/(P_1)_G\in E_{2'}$ by
induction on $|G|$ and thus $G\in E_{2'}$ by Lemma 3. So we may
assume that $(P_1)_G=1$. Since $P_1$ is $nc$-supplemented in $G$,
there is a subgroup $K$ of $G$ such that $P_1K\unlhd G$ and
$P_1\cap K\leq(P_1)_G=1$. Let $N=P_1K$. Then $G/N$ is solvable
since $|G/N|_2\leq 2$. We consider $N$. Since $P_1\cap K=1$ and
$P_1$ is a maximal subgroup of $P$, $|K|_2\leq 2$. Then $K$ has a
2-complement, $K_1$ say. Thus $K_1$ is also a 2-complement of $N$.
Then $N\in C_{2'}$ by \cite[Theorem A]{GC}. Hence $G$ has a Hall
$2'$-subgroup by Lemma 3 and thus $G\in E_{2'}$.

In the following $r$ is always a Mersenne prime. Since $G$ has a
Hall $2'$-subgroup, every composition factor of $G$ has a
2-complement. Hence every composition factor is either isomorphic
to $L_2(r)$ or a cyclic group of prime order by \cite[Corollary
5.6]{AF}.\pfend

{\bf Remark 2.} In (a) of Theorem 1, the hypothesis $P_1\neq 1$ is
necessary. For example, $G=L_2(7)$ and $p=7$. The example also
shows that the hypothesis $P_1\neq 1$ should be in the first
conclusion in \cite[Theorem 2.2]{FC}. Otherwise, From the example
it is certain that the conclusion that $[G:K]=p^r, r\geq 1$, in
the proof of \cite[Theorem 2.2]{FC}, is impossible.\pfend

{\bf Theorem 2.} Let $G$ be a finite group. Then $G$ is solvable
if and only if every Sylow subgroup of $G$ is $nc$-supplemented in
$G$.

\pf If $G$ is solvable, then every Sylow subgroup of $G$ has a
complement in $G$ and thus $nc$-supplemented in $G$.

Conversely, suppose that $G$ is a counterexample of smallest
order.

(1) If $N\unlhd G$, then $G/N$ is solvable.

Let $P$ be a Sylow $p$-subgroup of $G$, where $p$ is a prime
divisor of $|G|$. Then there is a subgroup $C$ of $G$ such that
$PC\unlhd G$ and $P\cap C\leq P_G$. Since $|PC:C|=|P:P\cap C|$ and
$P\in Syl_pG$, $(|PC:C|,|PC:P|)=1$. Hence $PC\cap N=(N\cap
P)(N\cap C)$ by \cite[A,1.2 Lemma]{DH}. Thus $NP\cap NC=N(P\cap
C)$ by \cite[A,1.6 Lemma (c)]{DH} and $PN/N\cap NC/N=N(P\cap
C)/N\leq P_GN/N\leq (PN/N)_{G/N}$. On the other hand,
$(PN/N)(CN/N)=(PC)N/N\unlhd G/N$. And for every Sylow $p$-subgroup
$S/N$ of $G/N$, we set $P\in Syl_pS$. Then $P\in Syl_pG$ and
$PN/N=S/N\in Syl_pG/N$. Therefore, $G/N$ satisfies the hypothesis
of the theorem. Then $G/N$ is solvable since $G$ is a
counterexample of smallest order.

(2) $G$ has a unique minimal normal subgroup $N$, $\Phi(G)=1$ and
$O_p(G)=1$ for any $p||G|$.

Since the class of all solvable groups is a saturated formation,
$G$ has only one minimal subgroup $N$ and $\Phi(G)=1$ by (1). If
$O_p(G)\neq 1$, then $G/O_p(G)$ is solvable by (1) and $G$ is
solvable, which contradicts that $G$ is a counterexample.

For any $p||G|$ and $P\in Syl_pG$, there exists a subgroup $C$ of
$G$ such that $PC\unlhd G$ and $P\cap C\leq P_G\leq O_p(G)=1$ by
our hypothesis and (2). Then $PC\geq N$ by (2) and $C$ is a
$p$-complement of $PC$. Thus $C\cap N$ is a $p$-complement of $N$
for any $p||G|$ and $N$ is solvable by \cite[I,3.5 Theorem]{DH}.
Hence $G$ is solvable, which contradicts that $G$ is a
counterexample. The final contradiction completes the proof.\pfend

{\bf Theorem 3.} Let $G$ be a finite group and let $P$ be a Sylow
2-subgroup of $G$. Suppose that every maximal subgroup of $P$ is
$nc$-supplemented in $G$. Then $G$ is solvable.

\pf Assume that $G$ is a counterexample of smallest order. In
particular, $G$ is non-solvable.

(1) $O_2(G)=1$ and $O_{2'}(G)=1$.

Assume that $O_2(G)\neq 1$. Then $G/O_2(G)$ either is of odd order
or satisfies the hypothesis of the theorem by Lemma 2. In the
first case $G/O_2(G)$ is solvable by the odd order theorem. In the
second case $G/O_2(G)$ is also solvable since $G$ is a
counterexample of smallest order. Thus in both cases $G/O_2(G)$ is
solvable, $G$ is also solvable, a contradiction.

Assume that $O_{2'}(G)\neq 1$. Then $G/O_{2'}(G)$ satisfies the
hypothesis of the theorem by Lemma 2 and thus $G/O_{2'}(G)$ is
solvable since $G$ is a counterexample of smallest order. In
addition $O_{2'}(G)$ is solvable by the odd order theorem again.
Hence $G$ is solvable, a contradiction.

(2) $G$ has a unique minimal normal subgroup $N$ and $N$ is a
direct product of some simple groups, which are isomorphic to
each other. Moreover, $G=PN$.

Let $N$ be a minimal normal subgroup of $G$. We consider $PN$.

We assume that $PN<G$. By Lemma 2 $PN$ satisfies the hypothesis of
the theorem, then $PN$ is solvable since $G$ is a counterexample
of smallest order. In particular, $N$ is solvable. Then either
$O_2(N)\neq 1$ or $O_{2'}(N)\neq 1$ and thus either $O_2(G)\geq
O_2(N)>1$ or $O_{2'}(G)\geq O_{2'}(N)>1$, which contradicts (1).
Now we in the case $PN=G.$ Then $G/N\cong P/P\cap N$ is solvable.
Since the class of all solvable groups is a saturated formation,
$G$ has a unique minimal normal subgroup $N$. Evidently $N$ is not
solvable and $N$ is a direct product of some simple groups, which
are isomorphic with each other.

(3) The final contradiction.

Let $P_1$ be a maximal subgroup of $P$. Then there is a $K\leq G$
such that $P_1K\unlhd G$ and $P_1\cap K\leq (P_1)_G\leq O_2(G)=1$
by (1). Thus $|K|_2\leq 2$, $K$ has a normal 2-complement
$K_{2'}$. Evidently, $K_{2'}$ is also a Hall $2'$-subgroup of
$P_1K$. In addition $N\leq P_1K$ by (2) since $P_1\neq 1$. Hence
$K_{2'}\cap N$ is a 2-complement of $N$. On the other hand $G=PN$.
So $K_{2'}$ is a 2-complement of $N$ and $G$. Set $N_2=P\cap N$,
$H=N_G(K_{2'})$ and $P'=P\cap H$. Then $N_2\in Syl_2N$,
$G=PN=PK_{2'}$ and $G=NH=N_2H$ by Frattini argument and
\cite[Theorem A]{GC}. So $P=P\cap(N_2H)=N_2(P\cap H)=N_2P'$. If
$G$ is simple, then $G$ is solvable by \cite[Corollary 3.2]{FC}
and Remark 1. Now we assume that $G$ is not simple. Then $P>N_2$
and $P'\neq 1$. Since $G=N_2H=PH$, $P'=P\cap H\in Syl_2H$. Since
$H<G$ by (1), $P>P'$. Then there is a maximal subgroup $P_1'$ of
$P$ such that $P'\leq P_1'$. Then there is a $K'\leq G$ such that
$P_1'K'\unlhd G$ and $P_1'\cap K'\leq (P_1')_G\leq O_2(G)=1$ by
(1). By the same argument as above, with $(P_1',K')$ in place of
$(P_1,K)$, we get: the normal 2-complement $K'_{2'}$ of $K'$ is
also a 2-complement of $N$ and $G$, and $P_1'K'\geq N$. Then
$P_1'K'=(P_1'N_2)K'=PK'=G$ since $P'\leq P_1'$. Since $G=PK_{2'}$,
we may assume $K_{2'}=K_{2'}'$ by \cite[VI,4.5]{Hu} and
\cite[Theorem A]{GC}. Then $K'\leq N_G(K_{2'}')=H$ and
$G=P_1'K'=P_1'H=P_1'P'K_{2'}=P_1'K_{2'}$. Hence
$|G|=|P_1||K_{2'}|<|P||K_{2'}|=|G|$, a contradiction. The final
contradiction completes the proof.\pfend

With Lemma 4, by the same argument as in the proof of theorem 3,
we get the following:

{\bf Theorem 4.} Let $G$ be a finite group and $P$ a Sylow
2-subgroup of $G$. If every 2-maximal subgroup of $P$ is
$nc$-supplemented in $G$ and $L_2(q)$-free,where $q\equiv 3(mod\
8)$ or $q\equiv 5(mod\ 8)$, then $G$ is solvable.\pfend

{\bf Remark 3.} With the condition \em $nc$-supplemented\em\ \ in
place of the condition \em $c$-supplemented\em\ \ in \cite[Theorem
3.1 and Theorem 4.2]{FC}, we can get that $G$ is solvable by
Theorem 3 and Theorem 4 respectively, but cannot conclude that
$G/O_p(G)$ is $p$-nilpotent. For example:

Let $G=H\wr\la a\ra$ is a wreath product of $H$ and $\la a\ra$,
where $H=Z_7\rtimes Z_3$ is a semi-direct product of $Z_7$ by
$Z_3$ but $H\ncong Z_7\times Z_3$, and $a=(1234567)$. Let $p=3$
and $P\in Syl_pG$. Then $P$ is an elementary Abelian $p$-subgroup.
For every subgroup $P_1$ of $P$ there exists a subgroup $C$ of $P$
such that $P_1\times C=P$. Let $F=O_7(G)$. Then
$F=Z_7\times\cdots\times Z_7$ (7 copies of $Z_7$). Hence $FC\leq
G$, $(FC)P_1=FP=H\lhd G$ and $FC\cap P_1=1$. Thus $P_1$ is
$nc$-supplemented in $G$. However, $O_p(G)=1$ and $G$ is not
$p$-nilpotent.\pfend

The following is related to $c$-supplemented subgroups.

{\bf Lemma A}\cite[Lemma 4.1]{FC}. Let $G$ be a finite group and
let $p$ be a prime divisor of $|G|$ such that $(|G|,p-1)=1$.
Assume that the order of $G$ is not divisible by $p^3$ and $G$ is
$A_4$-free. Then $G$ is $p$-nilpotent. In particular, if there
exists odd prime $p$ with $(|G|,p-1)=1$ and the order of $G$ is
not divisible by $p^3$. Then $G$ is $p$-nilpotent.\pfend

{\bf Remark 4}. From the hypotheses of Lemma A, we cannot get that
$G$ is $p$-nilpotent. For example:

Let $G=(Z_{19}\times Z_{19})\rtimes \la a\ra$ and $p=19$, where
$o(a)=5$ and $a\in Aut(Z_{19}\times Z_{19})$. Then $(G,p)$
satisfies the hypotheses of Lemma A. However, $G$ is not
$p$-nilpotent. From the example, it is certain that the conclusion
that $p=2$ and $q=3$, in the proof of \cite[Lemma 4.1]{FC}, is
impossible. But the mistake cannot affect the results and
arguments after \cite[Lemma 4.1]{FC} in \cite{FC}, since Lemma A
holds by \cite[Lemma 3.12]{GS} if $p$ is the smallest prime
divisor of $|G|$.\pfend

From the hypotheses of Lemma A, we get the following:

{\bf Theorem 5}. Let $G$ be a finite group and let $p$ be a prime
divisor of $|G|$ such that $(|G|,p-1)=1$. Assume that the order of
$G$ is not divisible by $p^3$ and is $A_4$-free. Then $G$ is
$p$-nilpotent or $G/O_{p'}(G)\cong (Z_p\times Z_p)\rtimes H$,
where $H\leq Aut(Z_p\times Z_p)$, and $H$ is a cyclic group whose
order is odd and divides $\frac{p+1}{2}$.

\pf If $O_{p'}(G)\neq 1$, then $G/O_{p'}(G)$ satisfies the
hypotheses and $G/O_{p'}(G)$ is $p$-nilpotent or $G/O_{p'}(G)\cong
(Z_p\times Z_p)\rtimes H$ by induction on $|G|$. So we assume
$O_{p'}(G)=1$.

If $|G|_p=p$, then $G$ is $p$-nilpotent by \cite[VI,2.6]{Hu},
since $$|N_G(P)/C_G(P)||(|Aut(P)|,|G|) =(p-1,|G|)=1,$$ where $P\in
Syl_p G$. Now we assume $|G|_p=p^2$. If $p$ is the smallest prime
divisor of $|G|$, then $G$ is $p$-nilpotent by \cite[Lemma
3.12]{GS}. So we may assume $p$ is odd. Then $|G|$ is odd since
$(|G|,p-1)=1$ and thus solvable by the odd order theorem. This
gives $O_{p}(G)\neq 1$ since $O_{p'}(G)=1$.

In case $|O_p(G)|=p$. Then $|G/O_p(G)|_p=p$ and $G/O_p(G)$ is
$p$-nilpotent as above. Let $T/O_p(G)$ is the normal
$p$-complement of $G/O_p(G)$. Then $|T|_p=p$ and thus $T$ is
$p$-nilpotent as above again. So $T$ has normal $p$-complement
$T_{p'}$ and $T_{p'}$ is a character subgroup of $T$. Then
$T_{p'}\unlhd G$ and $T_{p'}\leq O_{p'}(G)=1$. Thus $G$ is a
$p$-group and thus $p$-nilpotent.

In case $|O_p(G)|=p^2$. Then $O_p(G)\in Syl_p G$. Since
$O_{p'}(G)=1$, $F(G)=O_p(G)$ and then $C_G(O_p(G))=O_p(G)$
by\cite[A,10.6 Theorem]{DH}. Thus $G/O_p(G)\leq Aut(O_p(G))$. If
$O_p(G)$ is cyclic then $|Aut(O_p(G))|=p(p-1)$. Thus,
$|G/O_p(G)||(p-1,|G|)=1$ and $G=O_p(G)$ is a $p$-group. Now we
assume $O_p(G)$ is an elementary Abelian subgroup of $(p,p)$-type.
Then $|Aut(O_p(G))|=(p+1)p(p-1)^2$ and $|G/O_p(G)||\frac{p+1}{2}$
since $(|G|,p-1)=1$. In addition, $G$ has a subgroup $H$ such that
$G=HO_p(G)$ and $H\cap O_p(G)=1$ by Shur-Zassenhaus theorem.
Therefore, $G\cong (Z_p\times Z_p)\rtimes H$, where $H\leq
Aut(Z_p\times Z_p)$, $|H||\frac{p+1}{2}$ and $|H|$ is odd. In
addition, $|L_2(p)|=\frac{(p+1)p(p-1)}{2}$, $|Aut(Z_p\times
Z_p)|=|GL_2(p)|=(p+1)p(p-1)^2$, and $L_2(p)$ is a section of
$GL_2(p)$. Hence $H$ is isomorphic to a subgroup of $L_2(p)$ and
thus $H$ is cyclic by \cite[II,8.27]{Hu}.\pfend

\end{document}